\documentclass[11pt,noamsfonts]{amsart}

\usepackage{amsmath,amsfonts,amssymb,amscd,verbatim,xypic, xy}

\newcommand{\marginlabel}[1]%
  {\mbox{}\marginpar{\raggedleft\hspace{0pt}\bfseries\sf#1}}

\def\ZZ{{\mathbb Z}}

\def\QQ{{\mathbb Q}}
\def\PP{{\mathbb P}}

\def\cI{\mathcal{I}}

\def\cA{\mathcal{A}}

\def\cG{\mathcal{G}}

\def\cL{\mathcal{L}}
\def\cO{\mathcal{O}}

\def\cN{\mathcal{N}}

\def\cU{\mathcal{U}}

 \DeclareMathOperator{\Coker}{coker}
\DeclareMathOperator{\Ker}{ker} \DeclareMathOperator{\Pic}{Pic}

\DeclareMathOperator{\codim}{codim}

\newtheorem{lemma}{Lemma}[section]
\newtheorem{theorem}[lemma]{Theorem}
\newtheorem{corollary}[lemma]{Corollary}
\newtheorem{proposition}[lemma]{Proposition}

\theoremstyle{definition}
\newtheorem{definition}[lemma]{Definition}
\newtheorem{example}[lemma]{Example}

\newtheorem{remark}[lemma]{Remark}

\newtheorem{notation}{Notation}

\numberwithin{equation}{section}

\newcommand{\bean}{\begin{eqnarray}}
\newcommand{\eean}{\end{eqnarray}}
\newcommand{\be}{\begin{displaymath}}
\newcommand{\ee}{\end{displaymath}}
\newcommand{\bea}{\begin{eqnarray*}}
\newcommand{\eea}{\end{eqnarray*}}

\newcommand{\ol}{\overline}

\xyoption{arc}

\begin{document}

\title{Universal relations on stable map spaces in genus zero}

\author[ Anca M. Musta\c{t}\v{a}]{Anca~M.~Musta\c{t}\v{a}}
\author[Andrei Musta\c{t}\v{a}]{Andrei~Musta\c{t}\v{a}}
\address{Department of Mathematics, University of Illinois, 1409 W. Green Street, Urbana, IL 61801, USA} 
\email{{\tt amustata@math.uiuc.edu,
dmustata@math.uiuc.edu}}

\date{\today}

\begin{abstract}
We introduce a factorization for the map between moduli spaces of stable maps which forgets 
one marked point. This leads to a study of universal relations in the cohomology of stable
 map spaces in genus zero.
\end{abstract}
\maketitle

\bigskip

\section*{Introduction}

The moduli spaces of stable maps $\ol{M}_{0,m}(X,d)$ provide examples of Deligne-Mumford 
stacks whose intersection theory is both accessible and interesting, as indicated by the 
considerable success of Gromov-Witten theory in genus $0$. When the target is a point, $\ol{M}_{0,m}$
is a smooth projective variety and its cohomology ring has been computed by Keel (\cite{keel}). Recent studies lead to a more 
comprehensive view of the cohomology and Chow groups of these moduli spaces for other targets: \cite{behrend2}, 
\cite{getzlerpan}, \cite{oprea1}, \cite{cox1}, \cite{cox2}, \cite{noi1}, \cite{noi2}.  

 Let $d\in H_2(X)$ be a curve class on a smooth projective variety $X$. The space 
$\ol{M}_{0,0}(X,d)$ parametrizes maps from rational smooth or nodal curves into $X$ with image 
class $d$, such that any contracted component contains at least 3 nodes. Over $\ol{M}_{0,0}(X,d)$ 
there exists a tower of moduli spaces of stable maps with marked points and morphisms
$$f:\ol{M}_{0,m+1}(X,d)\to \ol{M}_{0,m}(X,d)$$ forgetting one marked point, such that $\ol{M}_{0,m+1}(X,d)$ 
is the universal family over $\ol{M}_{0,m}(X,d)$. 

In this paper we introduce a factorization of the forgetful map $f$, which gradually contracts 
part of the boundary. This allows a detailed study of the cohomology and Chow rings of $\ol{M}_{0,m+1}(X,d)$ as algebras over the rings of $\ol{M}_{0,m}(X,d)$. We find a series of universal relations over families of 
stable  maps. (Theorems 3.3 and 3.5 in text). 

This paper is the third in a series dedicated to the intersection rings of these moduli spaces. Previously we have found presentations for the Chow rings $\ol{M}_{0,m}(\PP^n,d)$ for all $m>0$ (\cite{noi1}, \cite{noi2}). The case $m=0$ seemed less accessible from our point of view, as in a sense  $\ol{M}_{0,m}(\PP^n,d)$ has more structure when $m>0$. A good parallel is in the study of the intersection ring for complete flag varieties versus that of the Grassmannian. In fact, when $d=1$, $\ol{M}_{0,0}(\PP^n,1)=Grass(\PP^1, \PP^n)$, while  $\ol{M}_{0,1}(\PP^n,1)$ is a flag variety.  This suggested an indirect approach, understanding $H^*(\ol{M}_{0,0}(\PP^n,d))$ by studying the extension of algebras  $H^*(\ol{M}_{0,0}(\PP^n,d))\to H^*(\ol{M}_{0,1}(\PP^n,d))$.  We view this step as prototypical for studies of Chow quotients by $SL_2$--action from an intersection--theoretical point of view. The description of $H^*(\ol{M}_{0,0}(\PP^n,d))$ is completed in \cite{noi3}. 

 The relation to Chow quotients in the sense of \cite{kapranov} is as follows. We may regard the coarse scheme $\ol{M}_{0,0}(X,d)$ as the $SL_2$-- Chow quotient of a simpler compactification for the space of maps $\PP^1\to X\times\PP^1$ with image class $(d,1)$: the space
$\Sigma_X^d$ of quasi-maps of $X$, also known as the "linear sigma model". This is true at least for X convex. Let $\alpha$ be the class of the generic $SL_2$--orbit in $\Sigma_X^d$. The Chow quotient of  $\Sigma_X^d$ is naturally a subvariety of the Chow variety $Chow ( \Sigma_X^d,\alpha)$. The graph space $\cU:=\ol{M}_{0,0}(X\times\PP^1,(d,1))$ is then the "universal family" over $\ol{M}_{0,0}(X,d)$ in the sense of the above, admitting a canonical map into $\Sigma_X^d$. Moreover, intermediate spaces between  $\cU$ and  $ \Sigma_X^d$, as described in \cite{noi1}, appear naturally when one considers the geometry of the $SL_2$--action. Finally, if $B$ be the Borel subgroup of $SL_2$, then $\ol{M}_{0,1}(X,d)$ is the quotient $\cU/B$. 
Our study may thus  be viewed as a procedure for understanding the intersection theory of $SL_2$ --Chow quotients $Y/ SL_2$ for a projective variety $Y$, under suitable assumptions on the stabilizers. The first step regards the relation between the cohomology of $Y$ and that of the quotient $\cU/B$, as in \cite{noi1}. The second step concerns the morphism  $\cU/B\to Y/ SL_2$, as in this paper. \cite{noi4} considers other applications of this viewpoint. 

 Another motivation for our study comes from the more complex situation of higher genus curves. Let $g>1$. Then by  methods similar to those presented in section 1 of this paper, the morphism 
  $$\ol{M}_{g,1}(\PP^n,d)\to  \ol{M}_{g,0}(\PP^n,d)\times_{\ol{M}_{g,0}}\ol{M}_{g,1}$$ 
factors in a series of birational morphisms with rational exceptional fibers. One may extract, for example, the Picard group for the basis in terms of that of the target. The line of inquiry opened in this paper partially applies to the maps above.  We note that some of the exceptional loci have codimension greater than 1, however, this situation is remedied when restricting over each of the strata in the natural stratification $\ol{M}_{g,0}=\bigcup M_{\gamma}$. Here $\gamma$ denote stable graphs and $M_{\gamma}$ are open sets of closed substrata $\ol{M}_{\gamma}$, representing specific splitting types of a genus $g$ curve. 
The case $g=1$ is special and requires separate treatment.

 The computations in this paper rely on two technical tools introduced in \cite{noi1}, \cite{noi2}. One tool is the intersection ring $B^*$ associated to a network of local regular embeddings, motivated by the structure of an \'etale cover of the target stack.
In the case of the moduli space of stable maps, the main advantage of $B^*$ over the usual intersection ring lies 
in the overall simplification of intersection on the boundary. Indeed, there is a natural stratification
of $\ol{M}_{0,m}(X,d)$, indexed by stable trees with degree and marking decorations. In contrast to 
$H^*(\ol{M}_{0,m}(X,d))$, $B^*$ allows the classes of each of these strata to be decomposed as a polynomial of divisor 
classes in a natural way.

The second tool is a set 
of stability conditions for maps to projective varieties which engender various spaces birational to $\ol{M}_{g,m}(X,d)$, via morphisms that contract the boundary of $\ol{M}_{g,m}(X,d)$. We call these the intermediate spaces of $\ol{M}_{g,m}(X,d)$. The fibered product of the universal families over these new spaces with $\ol{M}_{g,m}(X,d)$ interpolates between $\ol{M}_{g,m+1}(X,d)$
and $\ol{M}_{g,m}(X,d)$. 

 The sequence of spaces birational to $\ol{M}_{0,0}(\PP^n,d)$  as above is also constructed by Parker (\cite{parker}) via GIT quotients under the $SL_2$--action. Indeed, he regards $\ol{M}_{0,0}(\PP^n,d)$ as a GIT quotient of the graph space $\ol{M}_{0,0}(\PP^n\times\PP^1,(d,1))$ and its intermediate spaces as GIT quotients of the intermediate spaces for $\ol{M}_{0,0}(\PP^n\times\PP^1,(d,1))$ constructed in \cite{noi1}. The existence of this contraction of $\ol{M}_{0,0}(\PP^n,d)$ was  proved by different methods in \cite{chj}.

\section{A factorization of the forgetful map}

\subsection{} Let $m, n, d$ be natural numbers. We recall the various compactifications of $M_{0,m}(\PP^n,d)$ constructed in \cite{noi2}:

 Fix a rational number $a>0$ and an $m$-tuple $\cA=(a_1,...,a_m)\in \QQ^m$ such that $0\leq a_j\leq 1$ for all $j=1,...,m$ and such that $\sum_{i=1}^{m}a_i+da>2$.
\begin{definition}

 An $ (\cA, a)$-stable  family of degree $d$ nodal maps  from rational curves with $m$ marked sections to $\PP^n$  consists of  the following
 data:
\[( \pi \colon C \to S , \{p_i\}_{1\leq i\leq m}, \cL , e ) \]
where  $\cL$ is a line bundle on $C$ of degree $d$ on each fiber $C_s$, and
$e:\cO^{n+1}\to \cL$ is a morphism of fiber bundles (specified up to isomorphisms of the target) such that:

\begin{enumerate}
\item $\omega_{C|S}(\sum_{i=1}^ma_ip_i)\otimes \cL^a$ is relatively ample over $S$,
\item $\cG :=\Coker e$, restricted over each fiber $C_s$, is a skyscraper sheaf supported only on smooth points of $C_s$, and
\item for any $ p \in S$ and for any $ I \subseteq \{ 1,...,m \} $ (possibly empty) such that $ p= p_i$ for any $i \in S$ the following holds  $$\sum_{i \in I} a_i + a\dim\cG_{p}\leq 1.$$

\end{enumerate}

\end{definition}

\begin{proposition}(\cite{noi2})
The moduli problem of $ (\cA, a) $-stable degree $d$ nodal  maps
with $m$ marked points into $\PP^n$  is finely represented by a
smooth Deligne-Mumford stack $\ol{M}_{0,\cA}(\PP^n,d,a)$.

 For any other pair $ (\cA', a')$ as above such that $a_i\geq a_i'$ for all $i=1,...,m$ and $a\geq a'$, there is a natural birational morphism $$\ol{M}_{0,\cA}(\PP^n,d,a) \to \ol{M}_{0,\cA'}(\PP^n,d,a'),$$ a weighted blow-up along regular local embeddings.

\end{proposition}

As a corollary, we construct a sequence of spaces and morphisms
posed between $\ol{M}_{0,m+1}(\PP^n,d)$ and $\ol{M}_{0,m}(\PP^n,d)$. The case $m=0$ is somewhat special and we treat it separately.

\begin{corollary}
  There is a sequence of smooth Deligne-Mumford stacks $\{\ol{U}^k\}_{0\leq k\leq \lfloor  (d-1)/2 \rfloor }$, and of morphisms
$$\ol{M}_{0,1}(\PP^n,d)=\ol{U}^0\to...\to \ol{U}^{\lfloor (d-1)/2\rfloor }\to \ol{M}_{0,0}(\PP^n,d)$$
such that the morphism  $\ol{U}^{\lfloor (d-1)/2\rfloor }\to \ol{M}_{0,0}(\PP^n,d)$  has 1 dimensional fibers, and for each of the morphisms $f^{k-1}_k: \ol{U}^{k-1} \to \ol{U}^k$ there is one boundary divisor $\ol{U}^k_h$ of $\ol{U}^k$ with the property that $f^{k-1}_k$ restricted to  $\ol{U}^k_h$ has 1 dimensional fibers. 

    for each $k$ as above, there 

\end{corollary}

\begin{proof}
 Let $k$ be any natural number such that $0\leq k\leq \lfloor (d-1)/2 \rfloor$. Fix a rational number $\epsilon$ such that $0<\epsilon <1$. Consider the space $\ol{M}_{0,0}(\PP^n, d, 1/(k+\epsilon))$ constructed according to Definition 1.1, and its universal family $\tilde{U}^k$.  Alternatively, $\tilde{U}^k$ may be obtained directly by Definition 1.1, as the space $\ol{M}_{0,(0)}(\PP^n, d, 1/(k+\epsilon))$ of maps with one marked point of weight zero. Indeed, Definition 1.1 still makes sense for weights zero, although more than one marked point of zero weight may create singularities of the moduli space (\cite{hassett}).

By Proposition 1.2, the objects above are smooth Deligne-Mumford
stacks, and there are natural birational morphisms
          $$ \ol{M}_{0,0}(\PP^n, d, 1/(k'+\epsilon)) \to \ol{M}_{0,0}(\PP^n, d, 1/(k+\epsilon)) \mbox{ and } \tilde{U}^{k'}\to \tilde{U}^k$$
for $0\leq k'<k \leq \lfloor (d-1)/2 \rfloor$. Note that
$\ol{M}_{0,0}(\PP^n, d)=\ol{M}_{0,0}(\PP^n, d, 1/\epsilon), $ and
thus $\ol{M}_{0,1}(\PP^n, d) \cong \tilde{U}^0 $.

 The space $\ol{U}^k$ is defined as the fiber product $$\ol{U}^k=\ol{M}_{0,0}(\PP^n, d)\times_{\ol{M}_{0,0}(\PP^n, d, 1/(k+\epsilon))}\tilde{U}^k.$$ The existence of natural morphisms $\ol{U}^{k'}\to \ol{U}^k$ for $k'<k $ is implied by the universality property of fiber products.
The morphism from $\ol{M}_{0,1}(\PP^n, d)$ to $\ol{U}^k$ contracts
the components of degree no larger than $k$ in each fiber over
$\ol{M}_{0,0}(\PP^n,d)$.

\end{proof}

\begin{notation}
 For any $0\leq  k'\leq k
\leq \lfloor (d-1)/2 \rfloor$, we denote the natural morphisms constructed above by $f^{k'}_k
 : \ol{U}^{k'}\to \ol{U}^k$  and  $f^{k}:  \ol{U}^k \to \ol{M}_{0,0}(\PP^n,d)$.

\end{notation}

 Consider a homogeneous coordinate system $\bar{t}=(t_0:...:t_n)$ on $\PP^n$. An \'etale cover $\bigsqcup_{\bar{t} }\ol{M}_{0,0}(\PP^n,d,\bar{t})$ for the stack $\ol{M}_{0,0}(\PP^n,d)$ was constructed in \cite{fultonpan}.
$\ol{M}_{0,0}(\PP^n,d,\bar{t})$ parametrizes stable maps together
with $(n+1)d$ sections obtained by pullback of the hyperplanes
$(t_i=0)$ in $\PP^n$.

Pullback of $\ol{M}_{0,0}(\PP^n,d,\bar{t})$ to $\ol{U}^k$ yields
analogous \'etale covers $\bigsqcup_{\bar{t} }\ol{U}^k(\bar{t})$.
The morphisms $f^k_{k+1}:\ol{U}^k\to\ol{U}^{k+1}$ may be well
understood at the level these covers.

\begin{lemma}

  For $k<\lfloor (d-1)/2\rfloor $, $f^k_{k+1}:\ol{U}^k\to\ol{U}^{k+1}$ is a blow-up along a codimension 2 regular embedding.

 When $d$ is odd, the morphism $f^{\lfloor (d-1)/2\rfloor }:\ol{U}^{\lfloor (d-1)/2\rfloor} \to \ol{M}_{0,0}(\PP^n,d)$ is a $\PP^1$-- bundle.
 When $d$ is even, let $\ol{M}(d/2)$ denote the image of $\ol{M}_{0,1}(\PP^n, d/2)\times_{\PP^n}\ol{M}_{0,1}(\PP^n, d/2)$ in $\ol{M}_{0,0}(\PP^n,d)$.
Over the complement of $\ol{M}(d/2)$ in $\ol{M}_{0,0}(\PP^n,d)$,
$\ol{U}^{\lfloor (d-1)/2\rfloor }$ is a $\PP^1$--fibre bundle. Over
$\ol{M}(d/2)$, it is the $\ZZ_2$-- quotient of a stack obtained by gluing two
$\PP^1$-bundles along a common section.

\end{lemma}

\begin{proof}

Consider a cardinal $k$ a subset $h\subset \{1,...,d\}$. Let
$\bar{h}$ denote its complement and let
$\ol{M}_{h\bar{h}}:=\ol{M}_{0,1}(\PP^n,|h|)\times_{\PP^n}\ol{M}_{0,1}(\PP^n,|\bar{h}|)$.
We denote by $\ol{U}^k_{h\bar{h}}$ the fiber product
$\ol{M}_{h\bar{h}}\times_{\ol{M}_{0,0}(\PP^n,d)}\ol{U}^k$. It admits
a canonical section $s_h:  \ol{M}_{h\bar{h}}\to
\ol{U}^k_{h\bar{h}}$. Indeed, we recall from Definition 1.1 the existence of a bundle $\cL$ over $\ol{U}^k$, and of a morphism $e:\cO^{n+1}\to \cL$. Moreover, each the fiber $C_x$ of $f^k$ over a point $x\in\ol{M}_{h\bar{h}}$ contains a special point $p_x$ in the support of  $\Coker e_x$. Set theoretically, we may say $s_h(x):=p_x$. This extends canonically to a stack--theoretical definition.

A look at the \'etale covers establishes $\ol{U}^{k-1}$ as  the blow-up of $\ol{U}^{k}$
along the image of $s_h$. Moreover, while $\ol{U}^{k}_{h\bar{h}}\to \ol{U}^{k}$ is only a
 local regular embedding, the image of $s_h$ is regularly embedded in
$\ol{U}^{k+1}$.

When $d$ is odd, the existence of enough distinct sections makes the
morphism $f^{\lfloor (d-1)/2\rfloor }(\bar{t}):\ol{U}^{\lfloor
(d-1)/2\rfloor}(\bar{t}) \to \ol{M}_{0,0}(\PP^n,d,\bar{t})$ locally
trivial.   When $d$ is even, the same reasoning applies on the 
 complement of $\ol{M}(d/2)$.  Let now  $h\subset \{1,...,d\}$ be a cardinal $(d/2)$--subset, let $\ol{M}_{h\bar{h}}$ be defined  as above. Consider two copies  $\ol{M}_{h}$ and  $\ol{M}_{\bar{h}}$ of $ \ol{M}_{0,2}(\PP^n, d/2)\times_{\PP^n}\ol{M}_{0,1}(\PP^n, d/2)$. There is a fiber square diagram 
  \bea   \diagram \ol{M}_{h}\bigcup_s\ol{M}_{\bar{h}} \rto \dto &  \ol{M}_{0,1}(\PP^n,d)\times_{\ol{M}_{0,0}(\PP^n,d)} \ol{M}(d/2) \dto    \\
     \ol{M}_{h\bar{h}}    \rto   &    \ol{M}(d/2),     \enddiagram  \eea
where the two fiber products   $\ol{M}_{h}$ and  $\ol{M}_{\bar{h}}$  are glued along the canonical section $$s:\ol{M}_{0,1}(\PP^n, d/2)\times_{\PP^n}\ol{M}_{0,1}(\PP^n, d/2)\to  \ol{M}_{0,2}(\PP^n, d/2)\times_{\PP^n}\ol{M}_{0,1}(\PP^n, d/2).$$ There is a natural $\ZZ_2$--action on the two spaces at the left side of the diagram, and the quotients are on the right side. The same symmetry is preserved after $\ol{M}_{h}$ and $\ol{M}_{\bar{h}}$ are contracted to two $\PP^1$-- bundles over $\ol{M}_{h\bar{h}}$  by the successive blow-downs down to the $((d-2)/2)$--th step. (The precise moduli problems for these contraction are defined in section 2). 

\end{proof}
\bigskip

 \begin{notation} Let $A$, $B$ be two disjoint sets such that $A\bigcup B=\{1,...,m+1\}$ and let $i$ be a natural number, $0\leq i\leq d$. We denote by $D(A,i|B,d-i)$, the divisor representing split curves, of a degree $i$-- component containing the set $A$ of marked points, and
a degree $(d-i)$-- component containing the set $B$ of marked points.
\end{notation}

Similarly, there is a sequence of morphisms
$$\ol{M}_{0,m+1}(\PP^n, d)=\ol{U}^{0}_{m} \to  \ol{U}^{1}_{m}\to ...\to  \ol{U}^{d+m-2}_{m}\to  \ol{M}_{0,m}(\PP^n, d) $$
defined as follows. For $k=1,...,d+m-2$, let $a_k=1/(k+\epsilon)$, where $\epsilon$ is a small positive rational number, and let $\cA_k$ be the $m$-tuple consisting of one copy of $1$ and $m-1$ copies of $a_k$. The space $\ol{U}^{k}_{m}$ is the pullback to $\ol{M}_{0,m}(\PP^n, d)$ of the universal family over $\ol{M}_{0,\cA_k}(\PP^n, a_k)$. Let $1$ be the marked point of weight $1$ on the generic curve over $\ol{M}_{0,\cA_k}(\PP^n, a_k)$. Thus the morphism
$\ol{M}_{0,m+1}(\PP^n, d) \to \ol{U}^{1}_{m}$ is a blow-up with exceptional divisor $D(\{n+1\}, 1|\{1,...,n\}, d-1)$. When $k>1$,  the morphism $\ol{U}^{k-1}_{m} \to \ol{U}^{k}_{m}$ factors into a sequence of blow-ups along disjoint codimension two loci. Pullbacks of the exceptional divisors to $\ol{M}_{0,m+1}(\PP^n, d)$ are $D(A, i|B, d-i)$ such that $n+1\in A$, $1\in B$, and $|A|+i=k$. The corresponding blow-up locus is isomorphic to the support of $D(A\setminus\{n+1\}, i|B, d-i)$. The intermediate blow-ups  posed between   $\ol{U}^{k-1}_{m}$ and $\ol{U}^{k}_{m}$ are also moduli spaces in their own right. Indeed, such a space is obtained as pullback of the universal family over $\ol{M}_{0,\cA'_k}(\PP^n, a_k)$, where the  $m$-tuple $\cA'_k$ consists of one copy of $1$, some copies of $a_k$ and other copies of $a_{k-1}$, in the desired order. The universal family is in itself a moduli space of weighted pointed stable maps, where the generic point of its fiber is assigned weight $0$.

 Finally, $\ol{U}^{d+m-2}_{m}\to  \ol{M}_{0,m}(\PP^n, d)$ is a $\PP^1$-bundle. Indeed, let $C_x$ denote the fiber of the forgetful morphism $\ol{M}_{0,m+1}(\PP^n, d)\to \ol{M}_{0,m}(\PP^n, d)$ over the point $x$, with the $n$ marked sections $p_1, ...,p_n$. The fiber of $\ol{U}^{d+m-2}_{m}\to  \ol{M}_{0,m}(X, d)$ consists of the irreducible component of $C_x$ containing $p_1$. On rigid covers, the above morphism always has three disjoint sections: $p_1$, the universal section, and at least one hyperplane section.


\subsection{The case of a general target}

Let $X$ be a smooth complex projective variety, and $d\in H_2(X)$ a class curve. A factorization of the forgetful morphism $\ol{M}_{0,m+1}(X, d)\to \ol{M}_{0,m}(X, d)$ may be induced from any projective embedding of $X$. From the point of view of the boundary strata however, it is more natural to consider the embedding of $X$ in a product of projective spaces as follows. Consider $\cL_1, ..., \cL_s\in \Pic (X)$ very ample, such that their first Chern classes generate the algebraic part of $H^2(X,\ZZ)$. We consider the embedding of $X$ given by all $\cL_i$. To any curve class $d\in H_2(X)$ we assign an $s$-tuple $d=(d_1,...,d_s)$, such that $\int_d c_1(\cL_i)=d_i$.

Let $n_i:=h^0(\cL_i)-1$. Definition 1.1 above may be extended to the case when the target is $\prod_{i}{\PP^{n_i}}$ by considering an $s$-tuple of weights $a=(a^1,...,a^s)$, an $s$-tuple $\cL:=\{\cL'_i\}_i$ of line bundles with morphisms $e^i:\cO_C^{n_i+1} \to \cL'_i$, and defining $\cL^a:=\bigotimes_i\cL'^{\otimes a^i}_i$. The spaces $\ol{M}_{0,\cA}(\prod_{i}{\PP^{n_i}}, d, a)$ are still smooth Deligne-Mumford stacks, by the same reasons as for target $\PP^n$.

For each tuple of non-negative integers $k=(k_1,...,k_s, k'_2,...,k'_{m})$, define a system of weights $a_k:=(\frac{1}{k_1+\epsilon},...,\frac{1}{k_s+\epsilon})$ and $\cA_k:=(1, \frac{1}{k'_2+\epsilon},...,\frac{1}{k'_{m}+\epsilon}).$ To these we assign the moduli spaces  $\ol{U}^{k}_{m}(\prod_{i}{\PP^{n_i}})$ set between $\ol{M}_{0,m+1}(\prod_{i}\PP^{n_i}, d)$ and $\ol{M}_{0,m}(\prod_{i}\PP^{n_i}, d)$, constructed by the same method as above.
We define
$$\ol{U}^k_{m}(X):=\ol{U}^k_{m}(\prod_{i}{\PP^{n_i}})\times_{\ol{M}_{0,m}(\prod_{i}\PP^{n_i}, d)}\ol{M}_{0,m}(X, d).$$
 The morphism $\ol{M}_{0,m+1}(X,d)\to \ol{U}^k_{m}(X)$ contracts components $C$ of the fibers of $\ol{M}_{0,m+1}(X,d)\to\ol{M}_{0,m}(X,d)$ such that the degree of $C$ is $(l_1,...,l_s)$ and $\sum_{i=1}^s\frac{l_i}{k_i+\epsilon}+\sum_{p_j\in C}\frac{1}{k'_j+\epsilon}\leq 1$. Such a component will be called $k$--unstable.

\section{The extended cohomology and Chow rings}

 The factorization of the forgetful map facilitates understanding the structure of the ring $H^*(\ol{M}_{0,m+1}(\PP^n,d))$ as an algebra over $H^*(\ol{M}_{0,m}(\PP^n,d))$. Boundary classes are adjoined at each blow-up step, and codimension two relations among them exist. Additionally, $\ol{U}^{\lfloor (d-1)/2\rfloor}\to \ol{M}_{0,0}(\PP^n,d)$ admits a relative cotangent class  $\psi$', whose pullback is related to the canonical class $\psi$ on $\ol{M}_{0,1}(\PP^n,d)$.  The relation between intersection rings becomes truly transparent when we regard our spaces together with the networks generated by their canonical  closed strata. Then relations among cohomology classes may be deduced by simple induction.

Consider a general smooth projective target $X\hookrightarrow \prod_{i=1}^s\PP^{n_i}$ as in the preceding subsection, and $d=(d_1,...,d_s)$. We denote by $G=S_{d_1}\times...\times S_{d_s}$ the product of $s$ groups of permutations $S_{d_i}$.  We recall succinctly the construction of a $G$-network of morphisms  associated to $\ol{U}^{k}_m(X)$, and its extended cohomology ring. (see \cite{noi1} and \cite{noi2} for a more detailed motivation). We keep notations from the preceding section. Here we work in cohomology although the same construction works for Chow rings.

Let $I$ be a set whose elements are of the form $h=h_1\sqcup ...\sqcup h_s\sqcup M_h$ such that each $h_i\subset \{1,...,d_i\}$ and $M_h\subset\{2,...,m\}$.
 Assume that  $h\cap h'=h, $ $h'$ or $\emptyset$  for all $ h, h'\in I$. 


\begin{definition}
The space $\ol{U}^{k}_{m,I}(X)$
parametrizes degree $d$  stable maps $\varphi: (C, \{p_j\}_{j=1,...,m})\to X $
together with a smooth point $p_{m+1}\in C$ and marked closed curves $ \{C_{h}\}_{h\in I}$ such that:
\begin{enumerate}
\item   $p_{m+1} \not\in C'$ for any $k$--unstable curve $ C'\subset C$. 
\item $\forall h\in I$, $p_1\not\in C_h\subset C$, and the degree of the map $\varphi_{| C_{h}}$ is $(|h_1|,...,|h_s|)$.
\item the incidence relations among the elements of $I$ translate into analogous incidence relations among the curves $C_h$:
\begin{itemize}
\item $\forall h\in I$, $\forall i\in \{2,...,m\}$,  $p_i\in C_h$ iff $i \in M_h$.
\item $C_h\subset C_{h'}$ iff $h\subset h'$ and $C_h\bigcap C_{h'}=\emptyset$ iff $h\bigcap h'=\emptyset$.
\end{itemize}
 \end{enumerate}
 \end{definition}

The space  $\ol{U}^{k}_{m,I}(X)$ is locally embedded in  $\ol{U}^{k}_{m}(X)$. In fact, $\ol{U}^{k}_{m}(X)$ corresponds to the choice $I=\emptyset$. 
For any two subsets $I\subset J$ as above, there is natural local regular embedding $\phi^I_J: \ol{U}^k_{m,J}\to \ol{U}^k_{m,I}$. The spaces $\ol{U}^k_{m,I}$ with the morphisms $\phi^I_J$ form a network, on which the group $G$ acts naturally. 

 In the case of $\ol{M}_{0,0}(X, d)$ we employ slightly different notations. Here we take  $I$ to be a set of 2--partitions $h\bigcup \bar{h}= \bigsqcup_i \{1,...,d_i\}$ such that for any pair $(h,\bar{h}),$ $(h',\bar{h}') \in I$, the set $\{h\bigcap h', h\bigcap\bar{h}', \bar{h}\bigcap h', \bar{h}\bigcap \bar{h}'\}$ has exactly three elements. Then  $\ol{M}_{0,I}^k$ parametrizes  stable maps $\varphi: C\to X$, together with one point $p_1\in C$ and splittings $C_{h}\bigcup C_{\bar{h}}=  C$ for all $(h,\bar{h})\in I$ satisfying all the relevant properties from Definition 2.1.

 Let $\cI$ be the set of all sets $I$ as above. $G$ acts on $\cI$ by permutations. For each $I\in \cI$, let
$G_I\subset G$ be the subgroup which fixes all elements of $I$. Any
$g\in G$ induces a canonical isomorphism
$g:\ol{U}^k_{m,I}\to\ol{U}^k_{m,g(I)}$, such that the following  diagram is commutative
  \bea   \diagram      \ol{U}^k_{m,J}   \rto^g \dto^{\phi^I_J} &    \ol{U}^k_{m,g(J)}    \dto^{\phi^{g(I)}_{g(J)}} \\
         \ol{U}^k_{m,J}   \rto^g   &    \ol{U}^k_{m,g(J)}      \enddiagram  \eea
whenever $I\subset J$.


In this paper all cohomology is considered with rational coefficients. The extended cohomology ring $B^*(\ol{U}^k_m)$ is constructed as follows:

\begin{definition}

The graded $\QQ$-vector space $B^*(\ol{U}^k_m;\QQ)$ is
$$B^*(\ol{U}^k_m;\QQ):= \bigoplus_{l=0}^{\dim
(\ol{U}^k_m)}B^l(\ol{U}^k_m),$$ where the extended cohomology groups are \bea
B^l(\ol{U}^k_m)  :=\oplus_{I}
H^{l-\codim_{\ol{U}^k_m}\ol{U}^k_{m,I}}(\ol{U}^k_{m,I}) / \sim , \eea  the
sum taken after all subsets $I$ as above with $
\codim_{\ol{U}^k}\ol{U}^k_{m,I} \leq l$. The equivalence relation $\sim
$ is generated by: $$\phi^I_{J*}(\alpha)\sim \sum_{[g]\in
G_I/G_J} g_*(\alpha)$$ for any $\alpha \in H^*(\ol{U}^k_{m,J})$ and
$J\supset I$. Here $\phi^I_{J*}$ is the Gysin pushforward, as all $\phi^I_{J}: \ol{U}^k_{m,J}\to \ol{U}^k_{m,I}$ are local complete intersections. 
 Multiplication is defined as follows:
$$ \alpha\cdot_r\beta := \bar{\phi}^{I*}_{I\cup J}(\alpha) \cdot \bar{\phi}^{J*}_{I\cup
J}(\beta) \cdot c_{top}(\bar{\phi}^{I\cap J*}_{I\cup
J}\cN_{\ol{U}^k_{m,I\cap J}|\ol{U}^k_m})$$ for $\alpha \in
H^*(\ol{U}^k_{m,I})$ and $\beta \in H^*(\ol{U}^k_{m,J})$, where
$\cN_{\ol{U}^k_{m,I\cap J}|\ol{U}^k_m}$ denotes the normal bundle of
$\ol{U}^k_{m,I\cap J}$ in $\ol{U}^k_m$.
\end{definition}

The extended rings $B^*$ satisfy all the good properties of the usual intersection rings, under appropriate conditions. Thus  pullback is well defined for morphisms compatible with the $G$-network. An example is provided by the morphisms $f^{k'}_{k,m}: \ol{U}^{k'}_m\to \ol{U}^k_m$. For each $I$ as in Definition 2.1, consider the fiber square \bea   \diagram  \ol{U}^{k'}_{m}\times_{\ol{U}^k_{m}} \ol{U}^k_{m,I} \rto \dto^{f^{k'}_{k,m,I}} &  \ol{U}^{k'}_{m}  \dto^{f^{k'}_{k,m}}\\   \ol{U}^k_{m,I} \rto & \ol{U}^k_{m}\enddiagram  \eea
where  $\ol{U}^{k'}_{m}\times_{\ol{U}^k_{m}} \ol{U}^k_{m,I}$ is a union of normal strata in $\ol{U}^{k'}_m$.  The pullback $f^{k' *}_{k,m}: B^*( \ol{U}^k_m) \to B^*( \ol{U}^{k'}_m)$ is constructed by concatenating all pullbacks $f^{k' *}_{k,m,I}$. These are obviously compatible with the equivalence relation of Definition 2.2.

When the target $X$ is convex, all the moduli spaces $\ol{U}^k_{m,I}$ are smooth Deligne-Mumford stacks, and there is an isomorphism between  $B^*(\ol{U}^k_{m,I})$ and the analogously defined homology $B_*(\ol{U}^k_{m,I})$. Then pushforward is well defined for morphisms compatible with the $G$-network.

 As a central feature, $B^*(\ol{U}^k_m)$ contains the usual ring  $H^*(\ol{U}^k_m)$ as invariant subalgebra under the natural action of the group $G$. Working in 
 $B^*(\ol{M}_{0,m}(X,d))$  is very convenient because all boundary may be decomposed as polynomials of divisor classes. Indeed, for $h=h_1\sqcup ...\sqcup h_s\sqcup M_h$ as
 above, let the fundamental class of $\ol{U}^0_{m,h}$ in $B^*(\ol{M}_{0,m}(X,d))$ be denoted by $D_h$. Then  $[\ol{U}^0_{m,I}]=\prod_{h\in I}D_h$, while the image of $\ol{U}^0_{m,I}$ in $\ol{M}_{0,m}(X,d)$ has class $\sum_{g\in G}\prod_{h\in I}D_{g(h)}$. 

 Another important feature is that the product $D_h^2$ has a well known geometric significance. Indeed, the stratum of class $D_h$ is a fiber product of the form $\ol{U}^0_{m,h}=\ol{M}_h\times_{X}\ol{M}_{\bar{h}}$, where $\ol{M}_h$ and are themselves moduli spaces, $\ol{M}_{\bar{h}}$, and the product is along evaluation morphisms $ev_h:\ol{M}_h\to X$, $ ev_{\bar{h}}:\ol{M}_{\bar{h}}\to X$ There are canonical $\psi$-classes associated to  $ev_h$ and $ ev_{\bar{h}}$. Their pullbacks $\psi_h$ and $\psi_{\bar{h}}$ are known to satisfy the relation \bean \psi_h +\psi_{\bar{h}}=c_1(\cN^{\vee}_{\ol{U}^0_{m,h} | \ol{M}_{0,m}(X,d)}   )\eean on $\ol{U}^0_{m,h}$, which is $-D_h^2$ in $B^*(\ol{M}_{0,m}(X,d))$.

\section{Universal relations in cohomology and Chow rings}

We will denote  by $D_{h}$ the fundamental classes of boundary divisors in  $B^*(\ol{M}_{0,m}(X,d))$,  where $h=h_1\sqcup ...\sqcup h_s\sqcup M_h$ as above. In  $B^*(\ol{M}_{0,m}(X,d))$, the fundamental classes of boundary divisors are $D_{h\bar{h}}$, where $h=h_1\sqcup ...\sqcup h_s$. 

Consider the forgetful map $f:\ol{M}_{0,m+1}(\PP^n,d))\to \ol{M}_{0,m}(\PP^n,d))$ with $m$ canonical sections $s_i$, and the evaluation maps $ev_i: \ol{M}_{0,m+1}(\PP^n,d))\to \PP^n$, for $i=1,...,m+1$. Let $\psi_i:=c_1(L_i)$, where the tautological line bundle $L_i$ is pullback by $s_i$ of the relative cotangent bundle of $f$. For any class $\alpha\in H^*(X)$, we define the kappa class $k(\alpha):=f_*ev^*\alpha$. 

 $H$ will denote the hyperplane divisor on $\PP^n$. 

 We recall the following comparison formula on $\ol{M}_{0,m+1}(\PP^n,d)$
\bean  \psi_i=f_{m+1}^*\psi_i +D_{i,m+1},  \eean 
where $D_{i,m+1}$ is the Cartier divisor associated to the canonical section $\sigma_i$ of $f_{m+1}$ (see for example \cite{witten}). 

\begin{lemma}
The following codimension two class is zero in the cohomology ring
$H^*(\ol{M}_{0,1}(\PP^n,d))$: \bea &
R :=\psi^2-\sum_{h}N_{\psi h}  \psi D_h
+\sum_{(h,h')}N_{hh'}  D_hD_{h'}+\frac{3}{d^4} k(H^2)^2-\frac{4}{d^3}k(H^3), & \eea where
$$ N_{\psi h} =\frac{|h|^2}{d^2}(6-4\frac{|h|}{d}), $$
$$ N_{h'h} =N_{hh'} =\left\{ \begin{array}{l} \frac{|h|^2}{d^2}(6\frac{|h'|}{d}-2\frac{|h|}{d}-3\frac{|h'|^2}{d^2}) \mbox{ if
} h'\supseteq h, \mbox{ and } \\
 -3\frac{|h|^2|h'|^2}{d^4}\mbox{ if
} h'\cap h=\emptyset . \end{array} \right.$$

\end{lemma}

 Note that although the relation is written in $B^*(\ol{M}_{0,1}(\PP^n,d))$, it is in fact symmetrical in $\{h\}$, and thus it  is a relation in  $H^*(\ol{M}_{0,1}(\PP^n,d))$.
\begin{proof}
The formula can be checked by increasing induction on $d$. When
$d=1$,  $$\frac{\psi^2}{4}+\frac{3}{4}k(H^2)^2-k(H^3)=0$$ describes the usual
relation in  the cohomology of the flag space $F(*,\PP^1,\PP^n)$,
over that of the Grassmannian $G(\PP^1,\PP^n)$.

Assuming the formula true for all degrees less than $d$, we can
derive a codimension two relation on any divisor $D_h$. Indeed,
$D_h$ is the class of a normal stratum 
$\ol{U}_{h}:=\ol{M}_{0,1}(\PP^n,|h|)\times_{\PP^n}\ol{M}_{0,2}(\PP^n,|\bar{h}|)$,
where the fiber product is via evaluation maps at the first
marked points of $\ol{M}_{0,1}(\PP^n, |h|)$ and
$\ol{M}_{0,2}(\PP^n,|\bar{h}|)$. Let $\pi_h$  and $\pi_{\bar{h}}$ be
the projections from $\ol{U}_{h}$ on the first and second factors; let $f_h:
\ol{M}_{0,1}(\PP^n, |h|) \to \ol{M}_{0,0}(\PP^n,|h|)$, $f_{\bar{h}}:
\ol{M}_{0,2}(\PP^n,|\bar{h}|) \to \ol{M}_{0,1}(\PP^n,|\bar{h}|) $ be
the morphisms forgetting the first marked points. The evaluation maps at these points are $ev_h$ and $ev_{\bar{h}}$. The following $\psi$ classes on $\ol{U}_h$ will play a role in the computation:  $\psi_h:= \pi_h^*\psi_1$  and $\psi$, pullback of the $\psi$ class on $\ol{M}_{0,1}(\PP^n,d)$.

 We define $k_h(H^2):=\pi_h^*f_{h *}ev_h^*(H^2)$, $k_h(H^3):=\pi_h^*f_{h *}ev_h^*(H^3)$, $k_{\bar{h}}(H^2):=\pi_{\bar{h}}^*f_{\bar{h} *}ev_{\bar{h}}^*(H^2)$,
$k_{\bar{h}}(H^3):=\pi_{\bar{h}}^*f_{\bar{h} *}ev_{\bar{h}}^*(H^3)$.
 All the following arguments take place on $\ol{U}_{h}$. By the additivity of kappa classes (\cite{kk}, Lemma 3.3),
$$k_h(H^2)+k_{\bar{h}}(H^2)=k(H^2) \mbox{ and }
k_h(H^3)+k_{\bar{h}}(H^3)=k(H^3).$$
  The induction hypothesis gives a formula for $k_h(H^3)$ as a quadratic expression in $k_h(H^2)$, $\psi_h$ and divisors $\{D_{h'}\}_{h'\subset h}$. Simultaneously, pullback via $f_{\bar{h}}$ of the analogous relation on $\ol{M}_{0,1}(\PP^n,|\bar{h}|)$ expresses  $k_{\bar{h}}(H^3)$ as a quadratic function of $k_{\bar{h}}(H^2)$,  $\pi_{\bar{h}}^*f^{*}_{\bar{h}}\psi_2$,  and
$\pi_{\bar{h}}^*f^{*}_{\bar{h}}D_{h'}=D_{h'}+D_{h'\cup h},$ for sets
$h'$ such that $h'\cap h=\emptyset$. Furthermore, $\psi_h$ and $\pi_{\bar{h}}^*f^{*}_{\bar{h}}\psi_2$  may be written as:
$$\psi_h=\psi-\sum_{h'\supseteq h} D_{h'} \mbox{ and } \pi_{\bar{h}}^*f^{*}_{\bar{h}}\psi_2=\psi-D_{\bar{h}} $$
 (conform \cite{leepan}, \cite{noi1} Lemma 3.21, and a comparison formula (3.1)). Summation of $k_h(H^3)$ and $k_{\bar{h}}(H^3)$ yields $k(H^3)$ as an expression in $\psi$, $k_h(H^2)$, $k_{\bar{h}}(H^2)$, and boundary divisors.

  The classes $k_h(H^2)$ and $k_{\bar{h}}(H^2)$ can be written in terms of
$k(H^2)$, $\psi$ and boundary divisors via the divisorial relation in
Lemma 2.2 of \cite{pandharipande2}:  on $\ol{M}_{0,1}(\PP^n,d)$,
\bean \psi+\frac{2}{d}H-\frac{1}{d^2} k(H^2)-\sum_{h'}\frac{|h'|^2}{d^2}
D_{h'}=0. \eean

 Pullback by $\pi_{\bar{h}}\circ f_{\bar{h}}$ of the equivalent relation on
$\ol{M}_{0,1}(\PP^n,|\bar{h}|)$  is \bea
\psi-D_{\bar{h}}+\frac{2}{|\bar{h}|}H-\frac{1}{|\bar{h}|^2}k_{\bar{h}}(H^2)-\sum_{h'\subset
\bar{h}}\frac{|h'|^2}{|\bar{h}|^2} (D_{h'}+D_{h'\cup h})=0, \eea
 and thus, after eliminating $H$:  \bea
k_{\bar{h}}(H^2) =\frac{|\bar{h}|}{d}
k(H^2)-|\bar{h}||h|\psi-\sum_{h'\subseteq \bar{h} } \frac{|h'|^2|h|}{d}
D_{h'}-\sum_{h'\supseteq
h}|h|(\frac{|\bar{h}'|^2}{d}-|\bar{h}|)D_{h'}, \eea while
$k_h(H^2)=k(H^2)-k_{\bar{h}}(H^2)$.

 Substituting the expressions for $k_h(H^2)$ and $k_{\bar{h}}(H^2)$ into the above mentioned expression for $k(H^3)$,
 we obtain the desired formula modulo the codimension  two annihilator of $D_h$ in
$B^*(\ol{M}_{0,1}(\PP^n,d))$.
  By Theorem 3.23 in \cite{noi1}, this annihilator is generated by $ \psi
D_{\bar{h}}$, $D_{h'}D_{h''}$, $D_{h'}\psi$, $D_{h'}H$ where $h'$ is such that $h'\cap h\not= h', h \mbox{ or } \emptyset$, and
$\{D_{h'}(\psi-\sum_{h''\supseteq h\cup h'}D_{h''})\}_{h'\cap h
=\emptyset }$. Thus in $B^*(\ol{M}_{0,1}(\PP^n,d))$, $R =a(h)$,
where $a(h)$ is a linear combination of the terms above. Of
necessity, $a(h)=a(h')$ for any $h, h'\subset \{1,...,d\}$. But
there are no codimension two elements that annihilate all boundary
divisors in $B^*(\ol{M}_{0,1}(\PP^n,d))$. Indeed, the only
codimension two relations in $B^*(\ol{M}_{0,1}(\PP^n,d))$ are linear
combinations of monomials $D_{h'}D_{h''}$ for $h'\cap h''\not= h',
h'',\emptyset$. It is enough to consider $h$, $h'$ such that
$|h|=|h'|=d-1$. The only common annihilators for these are  of the
form $D_{h''}D_{\bar{h}''}$, where $h''\supset \bar{h}$ or
$h''\supset \bar{h}'$. But these do not annihilate $D_{h''}$.

\end{proof}

\begin{notation}

  Let $\ol{M}_{h\bar{h}}:= \ol{M}_{0,1}(\PP^n,|h|)\times_{\PP^n}\ol{M}_{0,1}(\PP^n,|\bar{h}|)$.  Consider the gluing map
$\ol{M}_{h\bar{h}}\to \ol{M}_{0,0}(\PP^n,d)$ and its class $D_{h\bar{h}}$ in
$B^*(\ol{M}_{0,0}(\PP^n,d))$. Let $\pi_h$ and $\pi_{\bar{h}}$ denote
the projections from $\ol{M}_{h\bar{h}}$ on the two components and let
$\psi_h:=\pi_h^*\psi_1$. The image of the class $\psi_h \in H^*(\ol{M}_{h\bar{h}})$ in
$B^*(\ol{M}_{0,0}(\PP^n,d))$ is a degree 2 class denoted by $F_h$. 
\end{notation}

 Based on the structure of $B^*(\ol{M}_{0,1}(\PP^n,k))$ found in \cite{noi1} Theorem 3.23,     it is not hard to see that $B^*(\ol{M}_{0,0}(\PP^n,d))$ is generated over $H^*(\ol{M}_{0,0}(\PP^n,d))$ by the set of classes $F_h$, $D_{h\bar{h}}$ for $h\subset\{1,...,d\}$. (In \cite{noi1} we have worked with Chow rings; the entire argument works identically for cohomology).

\begin{notation}
Choose $I:=\{h\subset \{1,...,d\}, |h|>d/2\}$ if $d$ is odd, and let $I$ additionally contain half of the sets $h$ with $|h|=d/2$ if $d$ is even, under the condition that no two sets $h$, $\bar{h}$ are simultaneously in $I$.
We define the following classes in $B^*(\ol{M}_{0,1}(\PP^n,d))$:
$$\psi'_I:= \psi-\sum_{h\in I}D_h, $$
 $$D_I(h):= \sum_{h' \in I, h'\subset h} D_{h'} \mbox{ and } D_I(\bar{h}):= \sum_{h' \in I, h'\subset h} D_{\bar{h}'}$$
 for any $h \in I$. The class $\psi_I (h)$ is defined as $\psi_I (h):=\psi'_I +D_I(h)$.

 When $d$ is odd, $\psi'_I$ is the pullback of the relative cotangent class for the morphism $\ol{U}^{\lfloor (d-1)/2\rfloor }\to \ol{M}_{0,0}(\PP^n,d)$.
\end{notation}

\begin{remark}
 The following relation between $\psi$ classes on $\ol{M}_{0,m}(\PP^n,d)$ by Y.P.Lee and R. Pandharipande (\cite{leepan}, Theorem 1) is instrumental in our computations $$\psi_i+\psi_j=D(i|j).$$
Here $m\geq 2$,  $i\mbox{ and } j\leq m$, and $D(i|j)$ is the divisor representing split curves, such that the marked points $i$ and $j$ lie in different components. With the notations from the previous sections, we will employ an analogous relation existing on intermediate spaces $\ol{U}^k_m$, where $m>0$ and $k=(k_1, \{a_i\}_{i=2,...,m})$ consists of a positive integer $k_1<d+m-1$, and weights  $a_i$ on the marked points. Indeed, the class $\psi_1$ on $\ol{M}_{m+1}(\PP^n,d)$ descends to $\psi'_1$ on each  $\ol{U}^k_m$, while for any $j=2,...,m+1$, the class $\psi_j-\sum_{j\in h} D_h$, descends to a class $\psi_j'$ on $\ol{U}^k_m$, where the sum above is after all divisors $D_h$ contracted by the morphism $\ol{M}_{m+1}(\PP^n,d)\to \ol{U}^k_m$, i.e. for all $h=h_1\sqcup M_h$ such that $|h|_k:=\frac{|h|}{k_1+\epsilon}+\sum_{i\in M_h}a_i\leq 1$. Thus on $\ol{U}^k_m$ the following relation holds
\bean \psi'_1+\psi'_j= \sum_{j\in h} D_h,\eean 
where the sum is taken after all $h$ with $|h|_k>1$.

\end{remark}

Next we show how the algebra $B^*(\ol{M}_{0,1}(\PP^n,d))$ may be constructed starting from $B^*(\ol{M}_{0,0}(\PP^n,d))$, by adjoining a divisor for each intermediate space $\ol{U}^{k}$. Each divisor comes with a natural quadratic relation. Thus, relation (1) below is pulled back from $\ol{U}^{\lfloor (d-1)/2\rfloor }$, while relations (2) and (3) are pulled back from  $B^*(\ol{U}^{|\bar{h}|-1})$. We keep notations from Lemma 3.1 throughout.

\begin{theorem}
  The algebra $B^*(\ol{M}_{0,1}(\PP^n,d))$  over $B^*(\ol{M}_{0,0}(\PP^n, d))$ is generated by divisors $\psi'_I$ and $\{D_{h}\}_{h\in I}$. The ideal of relations is generated by the following 
\begin{enumerate}
\item $\psi'^2_I-\sum_{h\not\in I}\frac{1}{2}N_{\psi h} (F_h-F_{\bar{h}}) +\sum_{h\cap h'=\emptyset }N_{hh'}
D_{h\bar{h}}D_{h'\bar{h}'}+\frac{3}{d^4} k(H^2)^2-\frac{4}{d^3}k(H^3)$ for a choice of the set $I$ as above.
\item $D_h^2-D_h(D_{h\bar{h}}-\psi_I(h))-\frac{1}{2} [ D_{h\bar{h}}(\psi_I(h)-D_I(\bar{h}))+F_h ]$ for all $h\in I$.
\item $D_h \ker \{B^*(\ol{M}_{0,0}(\PP^n,d))[\psi'_I, \{D_{h'}\}_{h'\in I; |h'|\leq |h|}]\to B^*(\ol{M}_{h\bar{h}})  \}$ where the above morphism sends $\psi'_I$  to $\psi_h+\sum_{h' \in I, h'\subset h}D_{h'\bar{h'}}$; it sends  $D_{h'}$ to $D_{h'\bar{h'}}$ if $h'\subset h$ and to $0$ otherwise.

\end{enumerate}

\end{theorem}

Note: We consider by convention $D_{\bar{h}h}=D_{h\bar{h}}$, such that both terms $D_{\bar{h}h}D_{h\bar{h}}$ and $D_{h\bar{h}}D_{\bar{h}h}$ appear in relation (1) above.

\begin{proof}

 \textbf{ Part I} We study $B^*(\ol{U}^{\lfloor (d-1)/2\rfloor })$ as an algebra over $B^*(\ol{M}_{0,0}(\PP^n,d))$.  When $d$ is odd, $\ol{U}^{\lfloor (d-1)/2\rfloor }$ is a projective bundle over $\ol{M}_{0,0}(\PP^n,d)$. Its cohomology ring is thus generated over $H^*(\ol{M}_{0,0}(\PP^n,d))$ by the cotangent class $\psi'$, which satisfies a degree 2 relation over the base ring. Pullback of $\psi'$ to $\ol{M}_{0,1}(\PP^n,d)$ is $\psi+\sum_{h\in I}D_h$. Substituting this in the relation of Lemma 3.1 yields (1).

 When $d$ is even, we refer to the description of $\ol{U}^{(d-2)/2 }$ from Lemma 1.4. 
Keeping the same notations, let $\ol{U}^{ (d-2)/2}(d/2):=\ol{U}^{ (d-2)/2}\times_{\ol{M}_{0,0}(\PP^n,d)}\ol{M}(d/2) $. We may apply the open--closed exact sequence on extended cohomology rings to $\ol{M}_{0,0}(\PP^n,d)=\ol{M}(d/2)\bigcup (\ol{M}_{0,0}(\PP^n,d)\setminus\ol{M}(d/2))$ and $\ol{U}^{ (d-2)/2 }=\ol{U}^{ (d-2)/2 }(d/2)\bigcup (\ol{U}^{ (d-2)/2}\setminus \ol{U}^{ (d-2)/2 }(d/2))$.  Over the complement of $\ol{M}(d/2)$, the cohomology ring of $\ol{U}^{ (d-2)/2 }\setminus \ol{U}^{ (d-2)/2 }(d/2) $ is generated by the relative cotangent class $\psi'=\psi-\sum_{|h|>d/2}D_h.$

 Let $h\subset\{1,...,d\}$ be a subset of cardinal $d/2$. Consider the fiber product $\ol{U}^{ |h|-1}_{h\bar{h}}:=\ol{U}^{|h|-1 }\times_{\ol{M}_{0,0}(\PP^n,d)}\ol{M}_{h\bar{h}}$. It is the union of two $\PP^1$ bundles  $\diagram \ol{U}^{ |h|-1}_{h}\rto^{f^{|h|-1}_h} & \ol{M}_{h\bar{h}} \enddiagram$ and $\diagram \ol{U}^{ |h|-1 }_{\bar{h}}\rto^{f^{|h|-1}_{\bar{h}}} & \ol{M}_{h\bar{h}} \enddiagram $ sharing a common section $S$. Let $s$ denote the class of $S$ in $B^*(\ol{U}^{|h|-1})$, as well as in $H^*( \ol{U}^{ |h|-1}_{h})$,  $H^*(\ol{U}^{ |h|-1 }_{\bar{h}})$. 

The classes of  $\ol{U}^{|h|-1}_{h}$ and $\ol{U}^{|h|-1 }_{\bar{h}}$ in
$B^*(\ol{U}^{|h|-1 })$  will be denoted by $D_{h}$ and
$D_{\bar{h}}$. Then $D_{h}+D_{\bar{h}}=D_{h\bar{h}}$ and $s=D_hD_{\bar{h}}$.  Thus, $\psi'$ and $\{ D_h\}_{|h|=d/2}$ are generators for the algebra $B^*(\ol{U}^{ |h|-1})$ over $B^*(\ol{M}_{0,0}(\PP^n,d))$.

 On each of the bundles  $\ol{U}^{ |h|-1}_{h}$ and $\ol{U}^{ |h|-1 }_{\bar{h}}$, the class $s$ may be written in terms of $\psi'$. 
Indeed,  $\ol{U}^{ |h|-1}_{h}= \ol{U}^{|h|-1}_1(\PP^n,d/2)\times_{\PP^n}\ol{M}_{0,1}(\PP^n,d/2),$ where  $\ol{U}^{|h|-1}_1(\PP^n,d/2)$ is a contraction of $\ol{M}_{0,2}(\PP^n,d/2)$ as described in subsection 1.2, such that the marked points both have weight 1, and the map to $\PP^n$ has weight $\frac{1}{d/2+\epsilon}$. The fiber product is along the evaluation maps at the first marked point. The $\psi_{\bar{h}}$ denote the $\psi$ class of  $\ol{U}^{|h|-1}_1(\PP^n,d/2)$ at the first marked point, as in Remark 3.2. 
 Then, by comparison formula (3.1) in conjunction with formula (3.3)  
 \bean -s=f^{|h|-1 *}_{h}\psi_{\bar{h}} - \psi_{\bar{h}}=f^{|h|-1 *}_{h}\psi_{\bar{h}} +
\psi'  \mbox{ on } \ol{U}^{ |h|-1}_{h}.  \eean  

  Moreover, as $f^{|h|-1}_{h\bar{h}}: \ol{U}^{|h|-1}_{h\bar{h}}\to \ol{M}_{h\bar{h}}$ is the sum of $f^{|h|-1 }_{h}$ and $f^{|h|-1 }_{\bar{h}}$ on components, then
$$ f^{|h|-1 *}_{h}\psi_{\bar{h}} =  f^{|h|-1 *}_{h\bar{h}}\psi_{\bar{h}}- f^{|h|-1
*}_{\bar{h}}\psi_{\bar{h}}$$ and, by formula (2.1),  \bean f^{|h|-1 *}_{h}\psi_{\bar{h}} =f^{|h|-1 *} F_{\bar{h}} +  f^{|h|-1 *}_{\bar{h}}\psi_{{h
}} + f^{|h|-1 *}_{\bar{h}}D_{h\bar{h}}.  \eean
 After comparing formulas (3.4) and (3.5) the their analogues for $\bar{h}$, we obtain the following relation in the algebra $B^*(\ol{U}^{ |h|-1})$ over $B^*(\ol{M}_{0,0}(\PP^n,d))$
  \bean F_h-F_{\bar{h}}= (D_{h}-D_{\bar{h}})(2 \psi'-D_{h\bar{h}}).\eean

This formula accounts for the different versions of relation (1) depending on the choice of the set $I$ above.  Indeed, choose $I$ such that $h\in I$. Then formula (3.6) is equivalent to 
 $$ F_h-F_{\bar{h}}= \psi'^2_{I\setminus\{h\}\cup\{\bar{h}\}}-\psi'^2_I,$$ which is exactly the difference between relation (1) applied to $I$ and to
$I\setminus \{h\}\bigcup \{\bar{h}\}$. Here we assumed the compatibility condition $D_hD_{h'\bar{h'}}=0$ for all $h, h'$ such that $|h|=|h'|=d/2$ (relation (3) in lemma). Note that in $B^*(\ol{M}_{0,0}(\PP^n,d))$, there is also a compatibility condition $D_{h\bar{h}}D_{h'\bar{h'}}=0$ for $h, h'$ as above.

 Equation (3.6) may be recast into relation (2) of the Lemma via formula (2.1). We note that the same equation may be obtained in part II of the proof by formally (at the level of \'etale covers) decomposing  $\ol{U}^{(d-2)/2}\to \ol{M}_{0,0}(\PP^n,d)$ into a blow-down along the section $S$ above, and a $\PP^1$-bundle. Moreover, the reasoning employed in part II of the proof also guarantees that there are no other relations besides (1), (2), and (3) (for $|h|=d/2$) in $B^*(\ol{U}^{(d-2)/2})$.

 \textbf{ Part II: The blow-down  $f^{k}_{k+1}:\ol{U}^{k}\to \ol{U}^{k+1}$. }
Once a choice of the set $I$ has been  fixed, we will drop the subscript $I$.

Consider a nonnegative integer $k < \lfloor (d-1)/2\rfloor$. We choose
$h\subset \{1,...,d\}$ such that $|h|=d-k-1$. With the notations
from Lemma 1.4,  consider the commutative diagram
\bea \xymatrix{ \ol{U}^{k}_{h} \ar[d]^{f^k_{k+1, h}}\ar@/^1pc/[rr]^{j^k_h} & \ol{U}^{k}_{\bar{h}} \ar[r]_{j^k_{\bar{h}}}\ar[d]^{f^k_{k+1, \bar{h}}} & \ol{U}^k \ar[d]^{f^k_{k+1}} \\ S_h\ar[r]^{j}\ar[dr]_{f^{k+1}_h} & \ol{U}^{k+1}_{h\bar{h}} \ar[r]^{j_{h\bar{h}}^{k+1}}\ar[d]^{f^{k+1}_{h\bar{h}}} & \ol{U}^{k+1} \ar[d]^{f^{k+1}}  \\ &  \ol{M}_{h\bar{h}}  \ar[r]^{j_{h\bar{h}}} & \ol{M}_{0,0}(\PP^n,d) }.\eea
$f_{h\bar{h}}^{k+1}$  admits a section $s_h:
\ol{M}_{h\bar{h}}\to \ol{U}^k_{h\bar{h}}$,  $S_h$ is the
image of $s_h$,  $j$ is the embedding of $S_h$ into
$\ol{U}^{k+1}_{h\bar{h}}$, and $j^{k+1}_{h}:= j_{h\bar{h}}^{k+1}\circ j$.
The  space $\ol{U}^{k}$ is  the blow-up of $\ol{U}^{k+1}$ along
$S_h$, and $\ol{U}^{k}_h$ is its exceptional
divisor, with regular embedding $j^k_h: \ol{U}^{k}_h\to
\ol{U}^{k}$. Let $f^k_h$
be the composition $  f^{k+1}_{h}\circ f^k_{k+1, h}$. In addition, $\ol{U}^{k}_{\bar{h}}$ is the strict transform of
$\ol{U}^{k+1}_{h\bar{h}}$ in $\ol{U}^{k}$ and $f^k_{\bar{h}}:= f^{k+1}_{h\bar{h}} \circ f^k_{k+1, \bar{h}}$.

 The fundamental class of $\ol{U}^{k+1}_{h\bar{h}}$ in  $B^*(\ol{U}^{k+1})$ is denoted by $D_{h\bar{h}}$, and is the pullback of the analogous class on $\ol{M}_{0,0}(X,d)$. The classes of $\ol{U}^{k}_h$ and $\ol{U}^{k}_{\bar{h}}$ in $\ol{U}^{k}$ are $D_h$ and $D_{\bar{h}}$, respectively.

 Let $\tilde{X}\to X$ be a blow-up along a regularly embedded locus $Y$. Following \cite{keel} (Theorem 2 in Appendix),  $H^*(\tilde{X})$ may be written explicitly as an algebra over $H^*(X)$, provided that the pullback morphism $H^*(X)\to H^*(Y)$ is surjective.
In the presence of compatible $G$--networks in  $X$, $Y$ and
$\tilde{X}$, the analogous statement holds for extended cohomology rings
(\cite{noi1}). Note that in the present case $H^*(\ol{U}^{k+1})\to H^*(S_h)$ is not surjective, while $j_h^{k+1 *}: B^*(\ol{U}^{k+1})\to B^*(S_h)$ is. 
We verify the surjectivity condition here.

 By induction, the algebra $B^*(U^{k+1})$ has generators $\psi'$ and $\{D_{h'}\}_{h'}$  over $B^*(\ol{M}_{0,0}(\PP^n,d))$, for all $h'\in I$ such that $|h'|<d-k-1$. On the other hand, the algebra $B^*(S_h)\cong B^*(\ol{M}_{h\bar{h}})$ over $B^*(\ol{M}_{0,0}(\PP^n,d))$ is generated by the divisor class $f_h^{k+1 *}\psi_h$. By formulas (3.7) and (3.8) below,  $f_h^{k+1 *}\psi_h=-j_h^{ k+1 *}\psi (h)$, while $j_h^{k+1 *}D_{h'}=D_{h'\bar{h'}}$ if $h'\subset h$ and $0$ otherwise. This proves that $j_h^{ k+1 *}$ is surjective.

By \cite{keel}, the ideal of relations in $B^*(\ol{U}^k)$ is made of two parts: $D_h\Ker j_h^{ k+1 *}$ and $D_h^2-aD_h+b$, where $a,b\in  B^*(\ol{U}^{k+1})$ are such that $j_h^{ k+1 *}a=c_1(\cN_{S_h|\ol{U}^{k+1}})$, and $b= j^{k+1}_{h *}[S_h]$. We proceed to find $a$ and $b$ in terms of the generators of $ B^*(\ol{U}^{k+1})$.

 Let $\psi(k+1)$ denote the first Chern class of the relative dualizing sheaf for $\ol{U}^{k+1}\to \ol{M}_{0,0}(\PP^n,d)$. Thus pullback of $\psi(k+1)$ to $\ol{M}_{0,1}(\PP^n,d)$ is $\psi'+\sum_{h'\in I, |h'|<|h|}T_{h'}$. 

The space $\ol{U}^{k+1}_{h\bar{h}}$ is a fiber product $\ol{U}^{k+1}_1(\PP^n,|h|)\times_{\PP^n}\ol{M}_{0,1}(\PP^n,|\bar{h}|)$, where $\ol{U}^{k'}_1(\PP^n,|h|)$ is a contraction of $\ol{M}_{0,2}(\PP^n,|h|)$ as described in 1.2, such that both marked points have weight 1 and the map to $\PP^n$ has weight $\frac{1}{k+\epsilon}$. The fiber product is taken along the first marked point.
  Let $\psi_h$ and $\tilde{\psi}$ on $\ol{U}^{k+1}_{h\bar{h}}$ be pullbacks of the two $\psi$-- classes from $\ol{U}^{k}_1(\PP^n,|h|)$.
  We note that $\tilde{\psi}$ differs from  $j_{h\bar{h}}^{k+1 *} \psi(k+1)$ by the class of the section $j_*[S_h]$.  Thus  by formula (3.3),     $$ \psi_h= - \psi(k+1)-j_*[S_h]+\sum_{h'\supset \bar{h}, |h'|<|h|}D_{h'}, $$ which may be reformulated as  \bean \psi_h= -\psi(h)+D(\bar{h})-j_*[S_h] \eean on   $\ol{U}^{k+1}_{h\bar{h}}$. On the other hand, by comparison formula (3.1) $$ -j_*[S_h]= f_{h\bar{h}}^{k+1 *}\psi_h-\psi_h.$$ 
 Putting these two equations together, we obtain
\bean j^{k+1}_{h *}[S_h]=-\sum_{|h|=d-k-1}
\frac{1}{2}[ F_h+ (\psi(h)-D(\bar{h}))D_{h\bar{h}} ]\eean
  Note that $j^*\psi_h=0$ on $S_h$. Thus the first Chern class $c_1(\cN_{S_h|\ol{U}^{k+1}})$ is $-(
\psi(h)-D_{h\bar{h}})),$ as $D(\bar{h})$ is also null on $S_h$. Relation
(2) follows.

\end{proof}

\subsection*{The case $m\geq 1$.} We have shown in section 1 how the forgetful morphism $f_{m+1}: \ol{M}_{0,m+1}(\PP^n,d)\to \ol{M}_{0,m}(\PP^n,d)$ factors out into a series of blow-downs and the projection of a $\PP^1$-- bundle over $\ol{M}_{0,m}(\PP^n,d)$. The first marked point was chosen to play a special role in our construction. The analysis done in the proof of Theorem 3.3. carries out to this case, with the extra simplification provided by the existence of sections in the intermediate spaces. Indeed, the $\PP^1$-- bundle over $\ol{M}_{0,m}(\PP^n,d)$ admits a canonical section $\sigma_1$. This determines a Cartier divisor $D_{1,m+1}$ on the bundle, whose class generates the cohomology of the bundle as an algebra over $H^*(\ol{M}_{0,m}(\PP^n,d))$. Note that  the subscript $\{1,m+1\}$ is a distinct convention from the subscripts $h$ employed elsewhere. The same notation will be used for the pull-back of the above divisor to $\ol{M}_{0,m+1}(\PP^n,d)$. The following well-known relation holds on $\ol{M}_{0,m}(\PP^n,d)$ (see for example \cite{pandharipande2}).
$$\sigma_1^*D_{1,m+1}=f_{m+1 *}(D_{1,m+1}^2)=-\psi_1,$$
while $\psi_1=f_{m+1}^*\psi_1+D_{1,m+1}$ on $\ol{M}_{0,m+1}(\PP^n,d)$.

 To every $h\subset\{1,...,d\}\bigsqcup \{2,...,m+1\}$ such that $m+1\in h$ and either $h\bigcap\{1,...,d\}\not=\emptyset$ or $|h\bigcap\{2,...,m+1\}|>2$, there corresponds a blow-down in the factorization of $f_{m+1}$.
 The blow-up locus is a section $S_h$ over $ \ol{M}_{h\setminus\{m+1\}}:=$
$$:=\ol{M}_{h\cap \{2,...,m\}, \bullet}(\PP^n,|h\bigcap \{1,...,d\}|)\times_{\PP^n}\ol{M}_{\{1,...,m\}\setminus h, \bullet}(\PP^n, |\{1,...,d\}\setminus h|),$$ whose class may be computed as above.
With the notations from section 1, the exceptional divisor $D_{h}$ is paired with a strict transform $D_{h\setminus\{m+1\}}$ of the divisor with the same name by the relation $f_{m+1}^*D_{h\setminus\{m+1\}}=D_{h}+D_{h\setminus\{m+1\}}$.
The following theorem follows by the same arguments as for Theorem 3.3.

\begin{theorem}
The algebra $B^*(\ol{M}_{0,m+1}(\PP^n,d))$  over $B^*(\ol{M}_{0,m}(\PP^n,d))$ is generated by the divisor classes $D_{1, m+1}$ and $\{D_{h}\}_{m+1\in h}$, where the sets $h \subset \{1,...,d\}\bigsqcup \{2,...,m+1\}$ are such that $h\bigcap\{1,...,d\}\not=\emptyset$ or $|h\bigcap\{2,...,m+1\}|>2$. The ideal of relations is generated by
\begin{enumerate}
\item $D_{1,m+1}^2+f_{m+1}^*\psi_1 \cdot D_{1,m+1}$;
\item $(D_h-f_{m+1}^*D_{h\setminus\{m+1\}})(D_h+f_{m+1}^*\psi_1+D_{1,m+1}-\sum_{\{m+1\}\cup h\subset h'}D_{h'})$, for $h$ as above;
\item $ D_h\ker (B^*(\ol{M}_{0,m}(\PP^n,d))\to B^*(\ol{M}_{h\setminus\{m+1\}})  )$, $D_hD_{1,m+1}$, and $D_hD_{h'}$ whenever $h\bigcap h'\not=h, h'$ or $\emptyset$.
\end{enumerate}
\end{theorem}




 The relations of Theorems 3.3 and 3.4 characterize in general the role of the forgetful map in the cohomology of $\ol{M}_{0,m}(X,d)$, for any smooth projective target $X$. Moreover, via a suitable embedding of $X$ in a product of projective spaces, these relations are refined to hold for all $h$ defined as in section 1.2. In this sense we call them universal. However, for an arbitrary convex target, the above theorems do not exhaust the list of generators for the algebras $B^*(\ol{M}_{0,m+1}(X,d))$ over $B^*(\ol{M}_{0,m}(X,d))$. Let us consider for example the case $m=0$. With the notations from Theorem 3.3, the issue is that in general, the morphism  $j_{h\bar{h}}^{k+1 *}: B^*(\ol{U}^{k+1}(X))\to B^*(\ol{M}_{h\bar{h}}(X))$ may not be surjective, for example, when the cohomology ring of $X$ is not generated by divisors. A class $\gamma_h\in B^*(\ol{M}_{h\bar{h}}(X))$ which is not in the image of  $j_{h\bar{h}}^{k+1 *}$ contributes the classes $f_{k+1, h}^{k *}\alpha_h$ and $f_{k+1, h}^{k *}\alpha_hD_h$ to $B^*(\ol{U}^{k}(X))$, while classes of the type $f_{k+1, h}^{k *}\alpha_hD_h^a$ with $a\geq 2$ can be written in terms of the above via relation (2) in Theorem 3.3.



Here we present a simple application to the case $X=\PP^n$, $d=2$.

\begin{example}
The ring $H^*(\ol{M}_{0,0}(\PP^n,2))$ was computed in
\cite{behrend2}.

By Theorem 3.3, $B^*(\ol{M}_{0,1}(\PP^n,2))$ is the algebra extensions of 
$H^*(\ol{M}_{0,0}(\PP^n,2))$  obtained by adjoining classes $D_1, D_2, \psi, f$ where 
$f:=\frac{1}{2}(F_1-F_2)$ and the following quadratic equations hold
 \bean (\psi-D_i)^2-\frac{3}{16} D^2 + \frac{3}{16} k(H^2)^2-\frac{1}{2}k(H^3)  \pm f, \eean
 \bean f^2+\frac{1}{2}D^4+\frac{3}{4} k(H^2)^2D^2- 2 k(H^3) D^2. \eean
Here $D=D_1+D_2$ is the class of $\ol{M}_{1,2}:=\ol{M}_{0,1}(\PP^n,1)\times_{\PP^n}\ol{M}_{0,1}(\PP^n,1)$, and the second equation is obtained by writing $F_1$ and $F_2$ in the two flag varieties $\ol{M}_{0,1}(\PP^n,1)$.

Thus $H^*(\ol{M}_{0,0}(\PP^n,2))$ may be recovered as the subring of
invariants of $B^*(\ol{M}_{0,1}(\PP^n,2))$ under an action
of $\ZZ_2\times\ZZ_2$. Extracting the invariant relations is a fun
exercise.

In \cite{noi1}, $B^*(\ol{M}_{0,1}(\PP^n,2))$ is computed as the
$\QQ$-algebra generated by divisors $H, \psi, D_1$ and $D_2$, with
the ideal of relations generated by:

$$H^{n+1} \mbox{ , }  D_1D_2 \psi \mbox{ , }   D_i(H +\psi )^{n+1} \mbox{ , } D_iQ(\psi-D_i) \mbox{ , }
\left.\sum_{i=1}^2 Q(s)\right|^{s=\psi-D_i}_{s=\psi}
+\frac{(H+2\psi)^{n+1}}{\psi},$$ where $$Q(t):=\frac{(H+\psi
+t)^{n+1}-(H+\psi)^{n+1}}{t}.$$

Let $D:=D_1+D_2$, $t:=\frac{1}{2}(k(H^2)-D)$ and
$b:=H+\psi=\frac{1}{4}(k(H^2)+D)$, such that
$Q(s)=\sum_{l=0}^{n}(b+s)^lb^{n-l}$. Let
$S_N:=\sum_{l=0}^{N}a_lb^{N-l}$, where
$a_l:=(b+\psi-D_1)^l+(b+\psi-D_2)^l-(b+\psi)^l+(b-\psi)^l$, such
that the last relation in the ring becomes $S_n=0$ and the first and
fourth imply $a_{n+1}=0$. Equivalently, $S_n=S_{n+1}=0$. A third
invariant relation is $Db^{n+1}=(t-2b)b^{n+1}=0$, obtained from the
third relation in $B^*(\ol{M}_{0,1}(\PP^n,2))$.

The following simple observation will come into play:

(*) If $a_1, ..., a_m$ are variables such that $\prod_{i=1}^m
a_i=0$, then for any polynomial $P$, the following relation holds:
$$\sum_{l=0}^m(-1)^l\sum_{i_1,..., i_l\in\{1,...,m\} }
P(\sum_{j\not\in\{ i_1,..., i_l\} }a_j) =0.$$

This helps us write $S_n$ invariantly. Indeed, by (*) applied to
$\psi, D_1$ and $D_2$,
$$a_l:=(b-D_1)^l+(b-D_2)^l+(b+\psi-D)^l+(b-\psi)^l-b^l-(b-D)^l.$$
 After applying recurrence relations for pairs of summands:
\bea & a_{l+2}-ta_{l+1}+(b^2-bD)a_{l}+ & \\ &
+D_1D_2[(b-D_1)^l+(b-D_2)^l]-(\psi-D)\psi
[(b+\psi-D)^l-(b-\psi)^l]=0 & \eea for all $l\geq 0$, which, via
relation $D_1D_2 \psi=0$  and observation $(*)$, yields $
a_{l+2}-ta_{l+1}+(b^2-bD-D_1D_2-(\psi-D)\psi)a_{l}=0$. Furthermore,
formulas (3.9) and (3.10) imply the following invariant formulation:
  $$k:=b^2-bD+D_1D_2-(\psi-D)\psi= \frac{1}{4}k(H^2)^2+\frac{1}{8}k(H^2)D+\frac{1}{8}D^2-\frac{1}{2}k(H^3),$$
and $a_{l+2}-ta_{l+1}+ka_{l}=0$, or, summing up:
$S_{l+2}-tS_{l+1}+kS_l=b^l(2b-t)$. The relations in
$H^*(\ol{M}_{0,0}(\PP^n,2))$ found above can be thus written as
$Y_{n+1}=0$, where
\bea  Y_{l+1}=\left( \begin{array}{l}b^{l+1}(2b-t)\\ S_{l+1}-tS_{l} \\ S_{l} \end{array} \right) = \left( \begin{array}{lll} b&0&0\\ 1&0&k\\
0&1&t \end{array} \right)^{l}\left( \begin{array}{l} b(2b-t) \\
2b-t\\2 \end{array} \right).\eea
 This is consistent with the ring structure computed in \cite{behrend2}.


\end{example}


\providecommand{\bysame}{\leavevmode\hbox
to3em{\hrulefill}\thinspace}

\end{document}